\documentclass[11pt, oneside]{article}
\usepackage{geometry} 
\usepackage{amsmath}
\usepackage[parfill]{parskip} 
\usepackage{graphicx}
\usepackage{amssymb}
\usepackage{amsthm}

\title{Groups Acting on Metric Spaces with Asymptotic Property C}
\author{Susan Beckhardt}
\date{}

\begin{document}
\maketitle

\begin{abstract}
We show that if a group $G$ acts by isometries on a metric space $M$ which has asymptotic property C, such that the quasi-stabilizers of a point $x \in M$ have asymptotic dimension less than or equal to $n$, then $G$ itself has asymptotic property C.
\end{abstract}

\textbf{Definition:} Let $\mathcal{F}$ be a collection of subsets of a metric space $X$. Say that $\mathcal{F}$ is \emph{uniformly bounded} if there exists $R>0$ such that for every $F \in \mathcal{F}$, $diam(F)<R$. For $r>0$, say that $\mathcal{F}$ is \emph{$r$-disjoint} if for any $F_1, F_2 \in \mathcal{F}$, $d(F_1,F_2) > r$.

\textbf{Definition:} Let $X$ be a metric space and $n \in \mathbb{N}$. We say that $X$ has \emph{asymptotic dimension at most $n$} (or $asdim(X)\leq n$) if for any $r>0$ there exist $n+1$ uniformly bounded families $\mathcal{F}_0, \mathcal{F}_1, \dots, \mathcal{F}_n$ of subsets of $X$ such that each family $\mathcal{F}_i$ is $r$-disjoint, and $\bigcup_{i=0}^n \mathcal{F}_i$ covers $X$.

We say that $X$ has \emph{finite asymptotic dimension} if there exists $n$ such that $asdim(X) \leq n$.

A weaker property of metric spaces is asymptotic property C:

\textbf{Definition:} A metric space $X$ has \emph{asymptotic property C} if for any sequence of real numbers  $0 < r_0 < r_1 < \dots$, there is some $n \in \mathbb{N}$ such that there exist uniformly bounded families $\mathcal{F}_0, \mathcal{F}_1, \dots, \mathcal{F}_n$ of subsets of $X$ such that each family $\mathcal{F}_i$ is $r_i$-disjoint, and $\bigcup_{i=0}^n \mathcal{F}_i$ covers $X$.

It follows immediately from the definition that any space with finite asymptotic dimension has asymptotic property C; however, the converse does not hold.

\textbf{Definition:} Let $\lambda > 0$ and let $X,Y$ be metric spaces. A function $f: X \to Y$ is \emph{$\lambda$-Lipschitz} if for all $x_1, x_2 \in X$, $d_Y(f(x_1), f(x_2)) \leq \lambda d_X(x_1, x_2)$.

\textbf{Definition:} Let $G$ be a group which acts by isometries on a metric space $M$. For $x \in M$ and $R > 0$, the \emph{$R$-quasi-stabilizer} of $x$ is $W_R(x) = \{ g \in G \mid d_M(x,gx) \leq R \}$.

For any such action by a group $G$ on a metric space $M$, for a fixed $x \in M$, there is a projection $\pi : G \to M$ given by $\pi (g) = g(x_0)$. Note that the $R$-quasi-stabilizer of $x_0$ is simply the inverse image of the (closed) ball of radius $R$ about $x_0$, i.e. $W_R(x_0) = \pi ^{-1} ( B_R(x_0))$. Also note that since the action of $G$ on $M$ is transitive, $\pi$ is surjective. 

Let $\lambda = max\{d(x_0, s(x_0)) \mid s \in S\}$; then $\pi$ is a $\lambda$-Lipschitz map. To see this, it suffices to show that for any $g \in G$ and $s \in S$, the inequality $d_X(\pi (gs), \pi(g)) \leq \lambda d_S(gs, g)$ holds. But since $G$ acts by isometries, $d_X(\pi (gs), \pi(g)) = d_X(gsx_0, gx_0) = d_X(sx_0, x_0)$, and $d_X(sx_0, x_0) \leq \lambda = \lambda d_S(gs, g)$ by the fact that the metric on $G$ is left-invariant and $s$ is a generator. Thus the desired inequality holds.

In \cite{bell-dran}, Theorem 2, Bell and Dranishnikov show that if a group acts by isometries on a metric space of finite asymptotic dimension, and all $R$-quasi-stabilizers have finite asymptotic dimension, then the group itself has finite asymptotic dimension. We prove a similar theorem in the more general case that the metric space has asymptotic property C.

\textbf{Proposition 1:} Let $G$ be a group with finite generating set $S$ which acts transitively by isometries on a metric space $M$, and fix a basepoint $x_0 \in M$. Suppose that $M$ has asymptotic property C and that for all $R > 0$, $W_R(x_0)$ has asymptotic dimension $\leq n$. Then $G$ has asymptotic property C.

\begin{proof}
Let $0 < r_0 < r_1 < r_2 < \dots$ be a sequence of real numbers.

Now since $M$ has asymptotic property C, we can choose families $\mathcal{F}_0, \mathcal{F}_1, \dots, \mathcal{F}_m$ of subsets of $M$ such that:
\begin{enumerate}
\item{For each $i$, $0 \leq i \leq m$, $\mathcal{F}_i$ is $\lambda r_{(i+1)(n+1)}$-disjoint,}
\item{All $\mathcal{F}_i$ are uniformly bounded by a number $R$,}
\item{$\bigcup_{i=1}^m \mathcal{F}_i$ is a cover of $M$.}
\end{enumerate}

Notice that by (1) and the fact that $\pi$ is $\lambda$-Lipschitz, each $\pi^{-1}(\mathcal{F}_i) = \{ \pi^{-1}(F) \mid F \in \mathcal{F}_i \}$ is $r_{(i+1)(n+1)}$-disjoint. Also, the families $\pi^{-1} (\mathcal{F}_i)$ together cover $G$, though in general their elements (subsets of $G$) are not bounded. 

We will have to further subdivide the elements of each family $\mathcal{F}_i$ using the $R$-quasi-stabilizer. Since asdim($W_R(x_0)$)$\leq n$, choose families $\mathcal{A}_0, \mathcal{A}_1, \dots, \mathcal{A}_n$ of subsets of $W_R(x_0)$ that are each $r_{(m+1)(n+1)}$-disjoint, uniformly bounded, and together cover $W_R(x_0)$.

For each $F \in \mathcal{F}_i$, $0 \leq i \leq m$, choose an element $g_F \in \pi^{-1}(F)$. Recalling that left multiplication is an isometry in $G$, notice that $g_F W_R(x_0)$ is isometric to $W_R(x_0)$. Furthermore, $\pi^{-1}(F) \subset \pi^{-1}(B_R(g_Fx_0)) = g_F W_R(x_0)$. This means that we can obtain a cover of $\pi^{-1}(F)$ by multiplying every element of the families $\mathcal{A}_0, \mathcal{A}_1, \dots, \mathcal{A}_n$ on the left by $g_F$. Doing the same for every $F$ yields a cover of $G$.

We break $G$ into $(m+1)(n+1)-1$ families of subsets $\mathcal{W}_k$ as follows: for every $0 \leq i \leq m$, $0 \leq j \leq n$, let \[\mathcal{W}_{i(n+1)+j} = \{ g_FA \cap \pi^{-1}(F) \mid F \in \mathcal{F}_i, A \in \mathcal{A}_j \}.\]

Notice that
\begin{enumerate}
\item{For each $k$, $0 \leq k \leq (m+1)(n+1)-1$, there is exactly one pair $(i,j)$ such that $k=i(n+1)+j$.}
\item{By the above argument, $\bigcup_{k} \mathcal{W}_k$ is a cover of $G$.}
\item{Each $\mathcal{W}_{i(n+1)+j}$ is uniformly bounded, since its elements are subsets of translations of elements of $\mathcal{A}_j$, which is uniformly bounded.}
\item{Each $\mathcal{W}_{i(n+1)+j}$ is $r_{i(n+1)+j}$-disjoint, by the following: Let \[g_FA \cap \pi^{-1}(F) \neq g_{F'}A' \cap \pi^{-1}(F') \in \mathcal{W}_{i(n+1)+j}.\] If $F \neq F'$, then $d_G(\pi^{-1}(F), \pi^{-1}(F') \geq r_{(i+1)(n+1)} \geq r_{i(n+1)+j}$, since $\pi^{-1}(\mathcal{F}_i$ is $r_{(i+1)(n+1)}$-disjoint. Otherwise, if $F = F'$ but $A \neq A'$, then $d_G(g_FA, g_FA') \geq r_{(m+1)(n+1)} \geq r_{i(n+1)+j}$, since $\mathcal{A}_j$ is $r_{(m+1)(n+1)}$-disjoint.}
\end{enumerate}

Therefore $G$ has asymptotic property C, as desired.
\end{proof}

A simple corollary uses a weaker property of metric spaces, straight finite decomposition complexity, first defined in \cite{dran-zar}.

\textbf{Definition:} Let $\mathcal{X}$ and $\mathcal{Y}$ be metric families (families of metric spaces), and let $r>0$. We say that $\mathcal{X}$ is \emph{$r$-decomposable over $\mathcal{Y}$} if for each $X \in \mathcal{X}$ there exist $r$-disjoint families $\mathcal{U}$ and $\mathcal{V}$ of subsets of $X$ such that $\mathcal{U} \cup \mathcal{V}$ is a cover for $X$ and $\mathcal{U} \cup \mathcal{V} \subseteq \mathcal{Y}$.

\textbf{Definition:} Let $\mathcal{X}$ be a family of metric spaces. We say that $\mathcal{X}$ has \emph{straight finite decomposition complexity (sFDC)} if for any sequence $0 \leq r_1 \leq r_2 \leq \dots$ of real numbers, there exists $n \in \mathbb{N}$ and metric families $\mathcal{X}=\mathcal{X}_0, \mathcal{X}_1, \dots, \mathcal{X}_n$ such that for each $1 \leq i \leq n$, $\mathcal{X}_i$ is $r_i$-decomposable over $\mathcal{X}_{i-1}$, and $\mathcal{X}_n$ is uniformly bounded.

In \cite{dran-zar}, Prop. 3.2, Dranishnikov and Zarichnyi show that asymptotic property C implies sFDC; hence the following corollary to Proposition 1 immediately follows:

\textbf{Corollary:} Let $G$ be a group with finite generating set $S$ which acts transitively by isometries on a metric space $M$, and fix a basepoint $x_0 \in M$. Suppose that $M$ has asymptotic property C and that for all $R > 0$, $W_R(x_0)$ has asymptotic dimension $\leq n$. Then G has straight finite decomposition complexity.

We believe this corollary should still hold under the hypothesis that $M$ has the weaker property of straight finite decomposition complexity rather than asymptotic property C. The details will be given in a future paper.

\textsc{Department of Mathematics and Statistics, SUNY, Albany, NY 12222}

\emph{E-mail address:} \texttt{sbeckhardt@albany.edu}

\end{document}